\documentclass[a4paper]{amsart}

\title{Assouad dimension of self-affine carpets}
\author{John M. Mackay}
\address{Department of Mathematics \\
 University of Illinois at Urbana-Champaign \\ Urbana IL.}
\curraddr{Mathematical Institute, 24-29 St Giles', Oxford OX1 3LB, UK.}
\email{john.mackay@maths.ox.ac.uk}
\date{\today}
\subjclass[2010]{Primary 28A78; Secondary 28A80, 37F35}
\keywords{Assouad dimension, conformal Assouad dimension, Bedford-McMullen carpets}

\usepackage{amsmath,amsthm,amsfonts,graphicx,amssymb}
\usepackage{hyperref} %[dvipdfm]

\numberwithin{equation}{section}
\newtheorem{theorem}[equation]{Theorem}
\newtheorem{proposition}[equation]{Proposition}

\newtheorem{lemma}[equation]{Lemma}

\DeclareMathOperator{\diam}{diam}

\newcommand{\eps}{\epsilon}

\newcommand{\cJ}{\mathcal{J}}
\newcommand{\cF}{\mathcal{F}}

\newcommand{\Cdim}{\mathcal{C}\mathrm{dim}}

\newcommand{\cM}{\mathcal{M}}

\newcommand{\ra}{\rightarrow}
\newcommand{\R}{\mathbb{R}}
\newcommand{\N}{\mathbb{N}}

\def\XXint#1#2#3{{\setbox0=\hbox{$#1{#2#3}{\int}$}
\vcenter{\hbox{$#2#3$}}\kern-.5\wd0}}

\numberwithin{equation}{section}

\begin{document}

\begin{abstract}
We calculate the Assouad dimension of the self-affine 
carpets of Bedford and McMullen, and of Lalley and Gatzouras.
We also calculate the conformal Assouad dimension 
of those carpets that are not self-similar.
\end{abstract}

\maketitle

%%%%%%%%%%%%%%%%%%%%%%%%%%%%%%%%%%%%%%%%%%%%%%%%%%%%%%%%%%%%%%%%%
\section{Introduction}\label{sec-intro}
Bedford and McMullen generalized the construction of the Sierpi\'nski
carpet to build a class of self-affine sets (``carpets'') in the plane
~\cite{Bed-84-Carpets,McM-84-Carpets}.
Their construction was later further generalized by
Lalley and Gatzouras~\cite{Lal-Gatz-92-self-affine}.
In this note we calculate the Assouad dimension of these carpets.
We also calculate their conformal Assouad dimension in the
non-self-similar case.  This calculation exhibits an interesting 
dichotomy: such a carpet is either minimal for conformal Assouad dimension,
or has conformal Assouad dimension zero.

We begin by considering the carpets of Bedford and McMullen.
Let us recall the construction of these sets.
Given integers $n \geq m$, and a fixed, non-empty set 
$A \subset \{0,1,\ldots,n-1\}\times\{0,1,\ldots,m-1\}$,
we can define the self-affine set 
\[
	S = S(A) = \left\{ \left( \sum_{i=1}^{\infty} \frac{x_i}{n^i}, 
		\sum_{i=1}^{\infty} \frac{y_i}{m^i} \right) : \forall i \in \N, (x_i,y_i) \in A \right\}.
\]
Following McMullen, we let $t_j$ be the number of elements $(i,j)$ of $A$, for each
row $0 \leq j < m$.  McMullen shows that the Hausdorff dimension of $S$ satisfies
\[
	\dim_H (S) = \log_m \Bigg( \sum_{j=0}^{m-1} t_j^{\log_n m} \Bigg).
\]
We let $s$ denote the number of rows which have an entry
in $A$, that is, $s = |\{j:t_j\neq 0\}|$.
McMullen demonstrates that the upper Minkowski dimension is given by
\[
	\overline{\dim}_M (S) = \log_m(s) + \log_n \left(\frac{|A|}{s}\right).
\]
(His result is stated for the upper Minkowski dimension, but his proof calculates
the Minkowski dimension as well.)

In the self-similar case ($n=m$) the carpet carries an Ahlfors regular
measure of dimension $\log_n(|A|)$, and
so the Hausdorff, upper Minkowski and Assouad
dimensions of the carpet all have this value.
It seems, however, that the Assouad dimension of $S$ has not been calculated
in the non-self-similar case ($n>m$).
(The definition of Assouad dimension is recalled in Section~\ref{sec-prelim}.)
\begin{theorem}\label{thm-main} When $n > m$, the Assouad dimension of $S$ is
	\[
		\dim_A ( S ) = \log_m(s) + \log_n(t),
	\]
	where $t = \max \{ t_j : 1 \leq j \leq m\}$.
\end{theorem}
Note that for a self-affine carpet with $n>m$, we have
\[
	\dim_H(S) < \dim_M(S) < \dim_A(S),
\]
unless we are in the ``uniform fibers'' case, that is, every non-zero $t_j$ equals $t$.

We prove Theorem~\ref{thm-main} in Section~\ref{sec-dima}.  The upper bound follows
from a straightforward counting argument.  To show the lower bound, we
build a suitable ``weak tangent'' to $S$ and use the scale-invariant properties
of Assouad dimension.

\begin{figure}
	\centering
	\fbox{\includegraphics[width=0.5\textwidth]{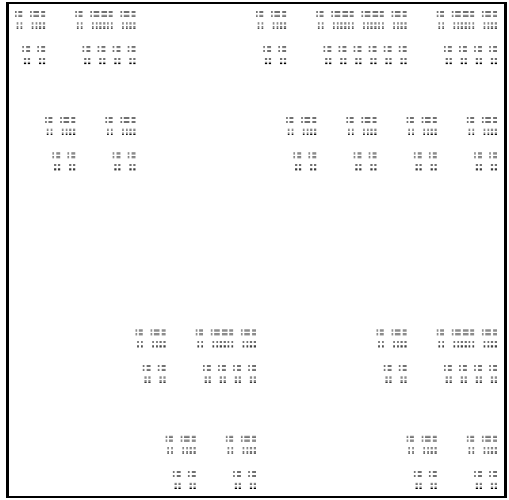}}
	\caption{Carpet $S_1$}\label{fig-43dust}
\end{figure}

\begin{figure}
	\centering
	\fbox{\includegraphics[width=0.5\textwidth]{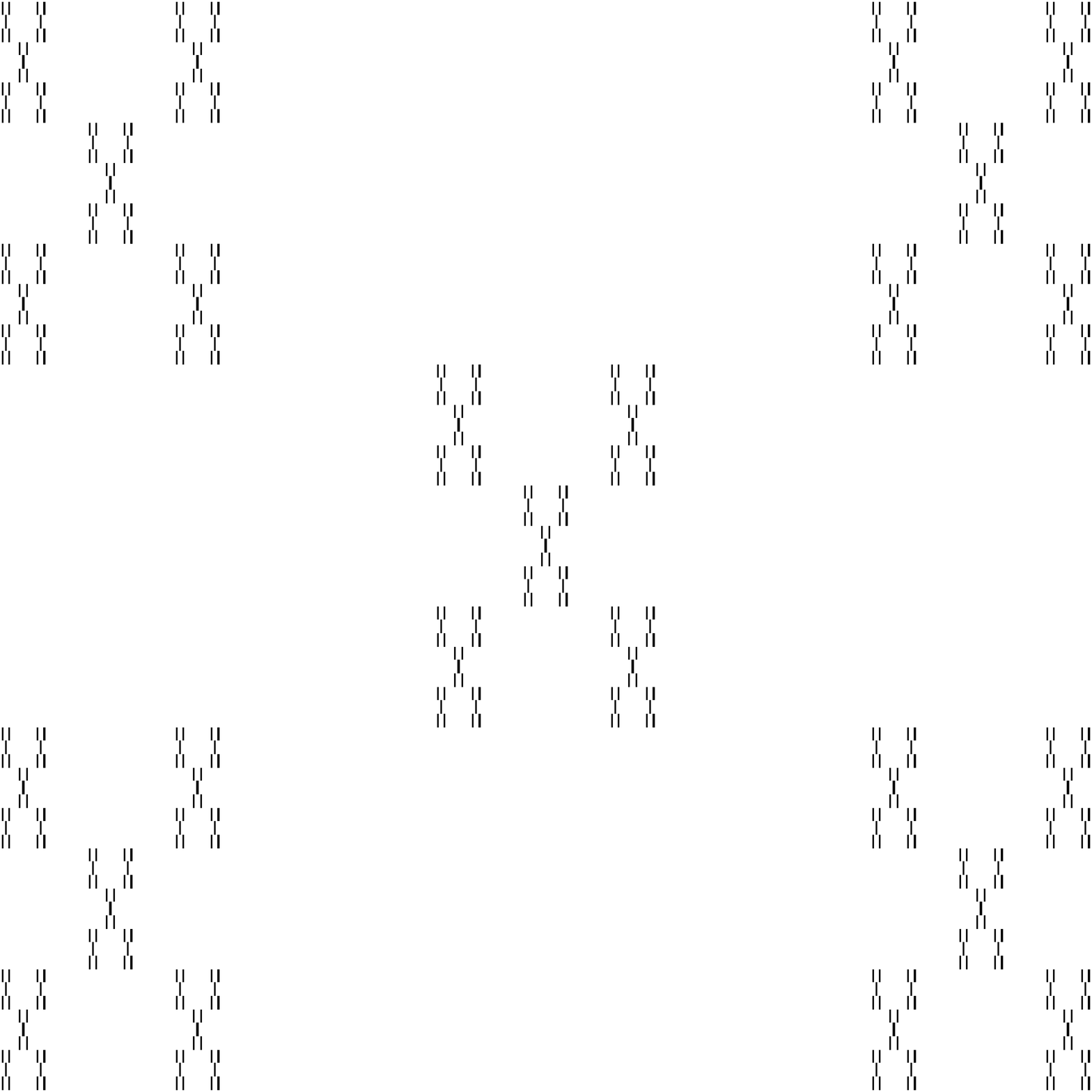}}
	\caption{Carpet $S_2$}\label{fig-53min}
\end{figure}

To illustrate this theorem, consider the carpet $S_1$ generated by
\[
	A_1 = \big\{ (0,2),(1,0),(2,2),(3,0),(3,2) \big\}, \ n=4,\ m=3,
\]
and the carpet $S_2$ generated by
\[
	A_2 = \big\{ (0,0),(0,2),(2,1),(4,0),(4,2) \big\},\ n=5,\ m=3.
\]
(See Figures~\ref{fig-43dust} and \ref{fig-53min} respectively.)
The theorem gives that $\dim_A(S_1) = \log_3(2)+\log_4(3)$ and
$\dim_A(S_2) = \log_3(3)+\log_5(2) = 1+\log_5(2)$.

Now, the Assouad dimension is a bi-Lipschitz invariant of a metric space,
but it may vary under quasi-symmetric deformations. (For example, quasi-conformal
homeomorphisms of the plane.)
The infimium of the values it can attain under these deformations is
called the conformal Assouad dimension of the metric space $X$,
and denoted by $\Cdim_A(X)$.  For more details
see \cite{Hei-01-lect-analysis,Mac-Tys-cdimexpo}.

Calculating the conformal dimension (Assouad or Hausdorff)
of a self-similar carpet is a challenging open problem 
(see, for example,~\cite{KL-04-confdim}).
In~\cite{Bin-Hak-cdim-carpets}, progress is made towards
calculating the conformal (Hausdorff) dimension
of self-affine ($n>m$) carpets.  
However, calculating the
conformal Assouad dimension of such carpets is quite simple.
\begin{theorem}\label{thm-cdima}
	Assume that $n>m$.  If both $t<n$ and $s<m$, then $\Cdim_A(S)=0$.
	Otherwise, $S$ is minimal for conformal Assouad dimension,
	i.e., $\Cdim_A(S) = \dim_A(S)$.
\end{theorem}

For our examples, we see that $\Cdim_A(S_1) = 0$, while
the carpet $S_2$ is minimal for conformal Assouad dimension.

The key observation in this result is that,
when $t=n$ or $s=m$, the weak tangent to $S$ built in the proof
of Theorem~\ref{thm-main} is the product of a Cantor set and an interval,
which is minimal for conformal Assouad dimension.
Since quasi-symmetric maps behave well with respect to taking tangents,
this gives the required bound.  See Section~\ref{sec-cdima} for details.

The methods and techniques of this paper apply to more general self-affine sets.
After the work of Bedford and McMullen, the Hausdorff and upper Minkowski dimension of 
more general sets were studied by
Lalley and Gatzouras~\cite{Lal-Gatz-92-self-affine}, Bara\'nski~\cite{Baranski-07-fractals}
and others.  For a recent survey on such constructions, see Chen and Pesin~\cite{Chen-Pesin-10-survey}.

In Section~\ref{sec-lal-gatz}, we extend Theorems~\ref{thm-main} and \ref{thm-cdima} to
the self-affine carpets of Lalley and Gatzouras.
Rather than specifying a collection of rectangles in a grid, as with the carpets of Bedford and McMullen,
the basic defining pattern of these carpets is a collection of $m$ disjoint 
rows of heights $b_1, b_2, \ldots, b_m$ in the unit square,
where the $i$th row contains $n_i$ disjoint self-affine copies of the entire set
of widths $a_{i1}, \ldots, a_{in_i}$.
We require that for every $1 \leq i \leq m$, $1 \leq j \leq n_i$, we have $a_{ij} < b_i$.
This pattern defines a self-affine set $S \subset [0,1]^2$.
(See Section~\ref{sec-lal-gatz} for more details.)

Let $\beta_y \in (0, 1]$ be the Hausdorff dimension of the projection of $S$ onto the $y$-axis,
namely, $\beta_y$ is the solution to $\sum_{i=1}^m b_i^{\beta_y}=1$.
Let $\beta_x \in (0, 1]$ be the maximal Hausdorff dimension of a horizontal fiber,
that is, $\beta_x = \max\{a : \exists i \text{ with } \sum_{j=1}^{n_i} a_{ij}^a = 1 \}$.
Then we show the following.
\begin{theorem}\label{thm-main-gl}
	The Assouad dimension of $S$ is $\dim_A(S) = \beta_x + \beta_y.$
\end{theorem}
\begin{theorem}\label{thm-cdima-gl}
	If $\beta_x<1$ and $\beta_y<1$, then $\Cdim_A(S) = 0$.
	Otherwise, $S$ is minimal for conformal Assouad dimension.
\end{theorem}

We expect that the results of this paper can be extended, at least partially, to even more general cases.
For example, the recent preprint of Bandt and K\"aen\-m\"aki~\cite{Ban-Kae-11-self-affine}
demonstrates that tangents to generic points in certain self-affine sets contain sets of the form
$C \times [0,1]$, where $C$ is a Cantor set.
It would be interesting to study how general a phenomenon the zero/minimal dichotomy 
for conformal Assouad dimension is amongst self-affine sets.

The author thanks Jeremy Tyson for introducing him to this topic and for
useful comments.  He also thanks Antti K\"aenm\"aki and the referee for their helpful suggestions.

\section{Preliminary results}\label{sec-prelim}
A metric space $X$ is \emph{doubling} if there exists an $N$ so that
any ball can be covered by $N$ balls of half the radius.  Repeatedly applying this
property, we see that there exists some $C>0$ and $\alpha>0$ so
that for any $r,R$ satisfying $0 < r \leq \frac{1}{2} R \leq \diam(X)$,
any ball $B(x,R) \subset X$ may be covered by $C (\frac{R}{r})^\alpha$ balls
of radius $r$.

The \emph{Assouad dimension} of a metric space $X$, denoted
by $\dim_A(X)$, is the infimal value of $\alpha$ for which there exists
a constant $C$ so that the above property holds.
We always have
\[
	\dim_H (X) \leq \overline{\dim}_M (X) \leq \dim_A (X),
\]
and these inequalities may be strict.
Unsurprisingly, if $X \subset \R^2$, then $\dim_A(X) \leq 2$.

Given $U \subset X$, the $\eps$-neighborhood of $U$ is the set
\[
	N(U,\eps) = \{ x \in X : \exists u \in U, d(x,u) < \eps\}
\]
Recall that the Hausdorff distance between $U,V \subset X$ is
the infimal $\epsilon$ so that $U \subset N(V, \eps)$ and
$V \subset N(U,\eps)$.  Denote this distance by $d_H(U,V)$.
If $X$ is compact, and $\cM(X)$ is the set of all closed subsets
of $X$, then $(\cM,d_H)$ is a compact metric space~\cite{BBI-01-metric-geom}.

We now use this convergence to give a non-trivial lower bound
on the Assouad dimension of a set.  For simplicity we restrict to
the case of subsets of $\R^2$, however this bound holds for general
``weak tangents''.
\begin{proposition}\label{prop-tangent-dimA}
	Fix a compact subset $X$ in $\R^2$.
	Suppose $U$ is a compact subset of $X$.
	Suppose that for each $k \in \N$, we have
	some $U_k \subset \R^2$ that is similar to $U$,
	i.e.\ $U_k$ is isometric to a possibly rescaled copy of $U$.
	Finally, suppose that $U_k \cap X$ converges to $\hat{U} \subset X$
	with respect to the Hausdorff distance.
	Then \[\dim_A (\hat{U}) \leq \dim_A(U).\]
	
	Moreover, $\Cdim_A(\hat{U}) \leq \Cdim_A(U)$.
\end{proposition}
\begin{proof}
	Suppose not.  Then there is some $\alpha$ so that
	$\dim_A(U) < \alpha < \dim_A(\hat{U})$.
	Then for all $D > 0$, there exists some $0 < r < R$ and
	a set $P$ in $\hat{U}$ of cardinality at least $D(\frac{R}{r})^\alpha$
	so that every pair of distinct points in $P$ 
	are separated by at least $r$.
	
	Since $U_k$ is similar to $U$, there is a fixed constant $C>0$
	so that every radius $R$ ball in $U_k$ can be covered by 
	$C( \frac{R}{r})^\alpha$ balls of radius $r$.
	On the other hand, for some sufficiently large 
	$k$ we can use $P$ to find a set 
	$Q \subset U_k \cap X \subset U_k$ that is 
	$\frac{r}{2}$-separated and lives in a ball of radius $2R$.  Therefore
	we require at least
	\[
		|Q| = |P| = D\left(\frac{R}{r}\right)^\alpha = 8^{-\alpha} D \left( \frac{2R}{r/4} \right)^\alpha
	\]
	balls of radius $\frac{r}{4}$ to cover $U_k$ inside this ball.  For sufficiently large
	$D$, this gives a contradiction.
	
	The lower bound for conformal Assouad dimension follows from the first part
	of the theorem and an Arzela-Ascoli type argument.  
	We sketch the argument for the reader's convenience; details are given
	in~\cite{Mac-Tys-cdimexpo}.
	
	Suppose that
	$\Cdim_A(U) < \Cdim_A(\hat{U})$.  Then there exists a quasi-symmetric homeomorphism
	$f:U \ra V$, where $\dim_A(V)<\Cdim_A(\hat{U})$.
	
	We can take a weak tangent to $f$ and get a quasi-symmetric map
	$\hat{f} : \hat{U} \ra \hat{V}$, where $\hat{V}$ is some weak tangent to $V$.
	Thus, by the first part of the theorem, we have a contradiction:
	\[
		\Cdim_A(\hat{U}) \leq \dim_A(\hat{V}) \leq \dim_A(V) < \Cdim_A(\hat{U}).\qedhere
	\]
\end{proof}

\section{Assouad dimension}\label{sec-dima}
\begin{proof}[Proof of Theorem~\ref{thm-main}]
First we prove the upper bound on the Assouad dimension.  This follows the
proof of the bound on upper Minkowski dimension given by McMullen~\cite{McM-84-Carpets}.

Since $n > m$, individual rectangles in the carpet 
get increasingly thin as we go down into the construction.  
To approximate squares with these rectangles we group 
them together as follows.
For any $k \in \N$, choose $l < k$ so that $n^l \leq m^k < n^{l+1}$.
That is, $l = \lfloor k\log_n(m) \rfloor$.  For any $p,q\in\N$, let
\begin{equation}\label{eq-Rk}
	R_k(p,q) = \left[ \frac{p}{n^l}, \frac{p+1}{n^l} \right] \times
		\left[ \frac{q}{m^k}, \frac{q+1}{m^k} \right].
\end{equation}

Let $\alpha = \log_m(s) + \log_n(t)$.  Since rectangles of the
form $R_k$ are present at every scale and location, and behave
like balls of radius $m^{-k}$, the proof that
$\dim_A(S) \leq \alpha$ reduces to the following lemma.
\begin{lemma}\label{lem-rects}
	There exists a constant $C$ so 
	that for every $1 \leq k' \leq k$, and any $p',q'$,
	the set $S \cap \mathrm{Int}(R_{k'}(p',q'))$ can be covered using
	at most $C m^{(k-k')\alpha}$ rectangles of the form $R_k(p,q)$.
\end{lemma}
\begin{proof}
	Let $l'$ and $l$ be chosen as before, corresponding
	to $k'$ and $k$ respectively.  Fix $R_{k'}(p',q')$.
	Let $N_k$ be the number of rectangles of the form
	$R_k(p,q)$ that meet $S \cap \mathrm{Int}(R_{k'}(p',q'))$.
	
	$N_k$ equals the number of ways to choose $(x_i)_{i=l'}^{l}$
	and $(y_i)_{i=k'}^{k}$, subject to certain restrictions.
	We have two cases to consider.
	
	\vspace{2mm}
	\noindent\textbf{Case 1:} $1 \leq l' \leq k' \leq l \leq k$.  Then
	\begin{enumerate}
		\item $(x_i, y_i) \in A$, $y_i$ is fixed by $q'$, for $i = l'+1, \ldots, k'$,
		\item $(x_i, y_i) \in A$, for $i = k'+1, \ldots, l$,
		\item $(\widetilde{x_i}, y_i) \in A$, for some $\widetilde{x_i}$, $i = l+1, \ldots k$.
	\end{enumerate}
	In this case
	\[
		N_k \leq (t)^{k'-l'} (|A|)^{l-k'} (s)^{k-l} 
			\leq t^{k'-l'} (st)^{l-k'} s^{k-l} = t^{l-l'} s^{k-k'},
	\]
	where we used that $|A| \leq st$.
	
	\vspace{2mm}
	\noindent\textbf{Case 2:} $1 \leq l' \leq l \leq k' \leq k$.  Then
	\begin{enumerate}
		\item $(x_i, y_i) \in A$, $y_i$ is fixed by $q'$, for $i = l'+1, \ldots, l$,
		\item $(\widetilde{x_i}, y_i) \in A$, for some $\widetilde{x_i}$, $y_i$ is fixed by $q'$,
			for $i = l+1, \ldots, k'$,
		\item $(\widetilde{x_i}, y_i) \in A$, for some $\widetilde{x_i}$, $i = k'+1, \ldots k$.
	\end{enumerate}
	Again we see that
	\[
		N_k \leq (t)^{l-l'} (1)^{k'-l} (s)^{k-k'} = t^{l-l'} s^{k-k'}.
	\]
	
	Therefore,
	\begin{align*}
		\log_m (N_k) & \leq \log_m( t^{l-l'} s^{k-k'} ) = (l-l') \log_m(t) + (k-k')\log_m(s) \\
			& \leq \Big(k\log_n(m) - k'\log_n(m)+1\Big) \log_m(t) + (k-k')\log_m(s) \\
			& = (k-k')\Big(\log_n(t) + \log_m(s)\Big) + \log_m(t).\qedhere
	\end{align*}
\end{proof}

It remains to bound the Assouad dimension from below.  
Choose $y_*$ with $0 \leq y_* < m$ so that $t = t_{y_*}$.
Fix some $x_*$ so that $(x_*, y_*) \in A$.
We will follow this rectangle into the construction 
in order to build a suitable weak tangent.

For each $k \in \N$,
let $p_k = \sum_{i=0}^{l-1} x_* n^i$, and
let $q_k = \sum_{i=0}^{k-1} y_* m^i$, where $l$ is
related to $k$ as before.
Let $f_k : \R^2 \rightarrow \R^2$ be the similarity
with scaling factor $m^k$ that takes $R_k(p_k,q_k)$
to $[0,m^k n^{-l}] \times [0,1]$.
Note that $m^k n^{-l} \in [1,n)$.

Fix $X = [0,n+1] \times [0,1]$.  Since $(\cM(X),d_H)$
is compact, some subsequence of $f_k(S) \cap X$ converges
to a compact set $\hat{S}$ in $X$.  Furthermore, we can assume that
$m^k n^{-l}$ converges to some $w \in [1,n]$.
Let $g:\R^2\ra \R^2$ be given by $g(x,y) = (x/w, y)$, and let 
\[
	W = g\big( \hat{S} \cap ([0,w] \times [0,1])\big) \subset [0,1]^2.
\]

$W$ looks like the product of two Cantor sets.  To be precise,
define $A'$ to be the set of pairs $(x,y)$ which satisfy
$(x, y_*) \in A$ and $(\tilde{x},y) \in A$, for some $\tilde{x}$,
and for $y_*$ as fixed above.

We can use $A'$ to build another carpet $S'=S(A')$.  This carpet has
the same $s$ value as before, and each non-empty row has $t$ entries.
Therefore McMullen's result shows that
\begin{equation}\label{eq-dimM-tangentcarpet}
	\overline{\dim}_M(S') = \log_m(s) + \log_n(t).
\end{equation}
In fact, by construction $S'$ is the product of two (self-similar) Cantor 
sets $C_x$ and $C_y$,
of dimensions $\log_n(t)$ and $\log_m(s)$ respectively.
See Figure~\ref{fig-53zoom} for an enlarged part of carpet $S_2$ from the introduction,
showing part of this structure emerging.

\begin{figure}
	\centering
	\fbox{\includegraphics[width=0.5\textwidth]{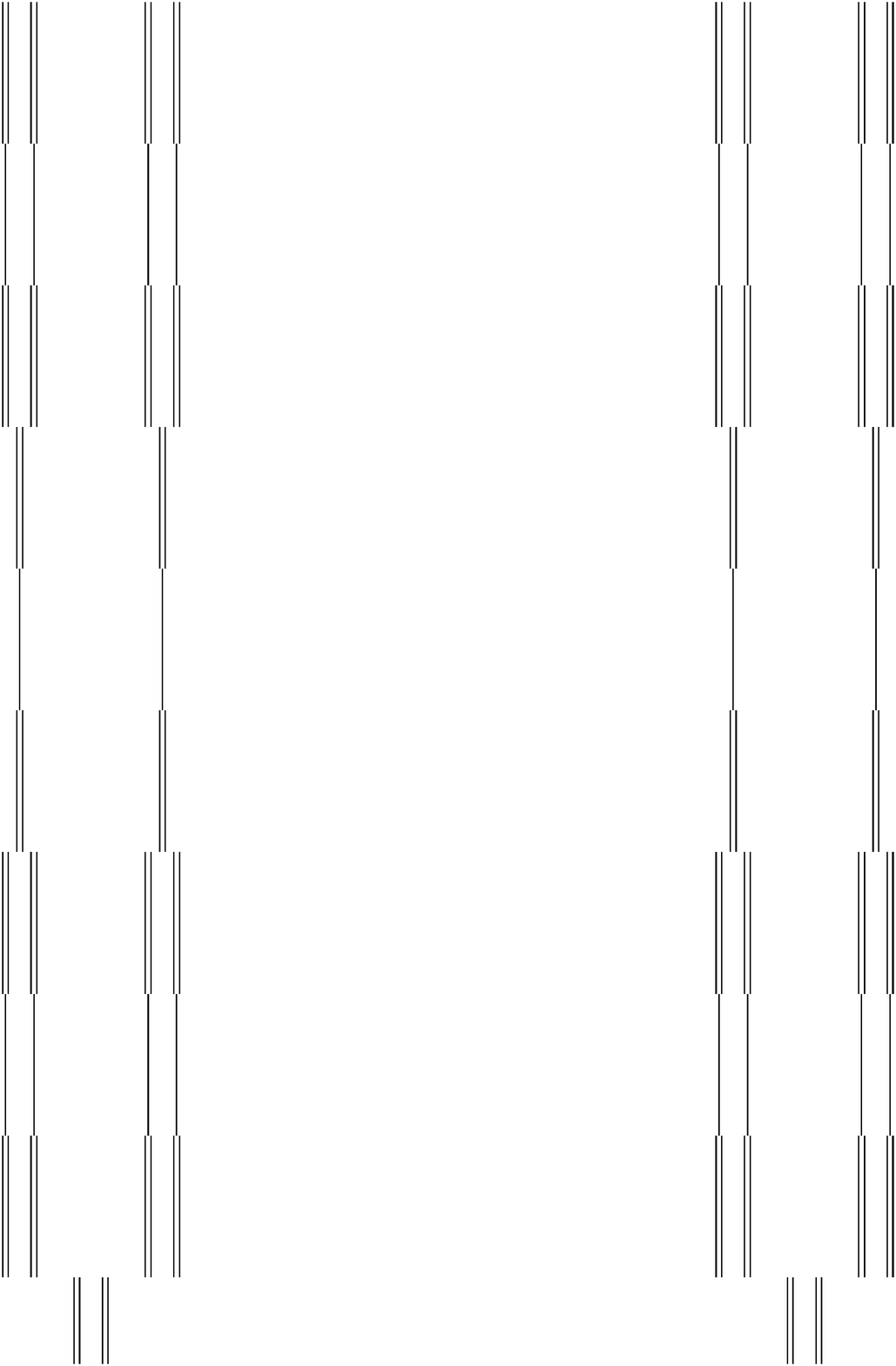}}
	\caption{A magnified part of Carpet $S_2$}\label{fig-53zoom}
\end{figure}

\begin{lemma}\label{lem-carpet-tangent}
	With the above notation, $W = S'$.
\end{lemma}
  
\begin{proof}
	Consider rectangles in $[0,1]^2$ of the form
	\begin{equation}\label{eq-special-rect}
		\left[ \sum_{i=1}^{k-l} \frac{x_i}{n^i} , 
			\sum_{i=1}^{k-l} \frac{x_i}{n^i} + n^{-(k-l)} \right] \times [0,1],
	\end{equation}
	where $(x_i, y_*) \in A$ for each $1 \leq i \leq k-l$.
	Let $W_k$ be the subset of $[0,1]^2$ given by placing an affine copy
	of $S$ into each such rectangle.
	Note that $W_k$ is just an affine copy of the set 
	$R_k(p_k, q_k) \cap S$ (with uniformly bounded distortion).
	Consequently, $W$ is the Hausdorff limit of the sets $W_k$.
	
	Now, each rectangle in \eqref{eq-special-rect} is of width $n^{-(k-l)}$,
	so the copy of $S$ inside is within Hausdorff distance $n^{-(k-l)}$ of a
	copy of $C_y$ given the appropriate $x$-coordinate.  
	Moreover, the $x$-coordinates of the rectangles in \eqref{eq-special-rect}
	are within Hausdorff distance $n^{-(k-l)}$ of $C_x$.
	Therefore,
	\[ d_H(W_k, C_x \times C_y)  \leq 2 n^{-(k-l)} 
		\leq 2 \left(\frac{m}{n}\right)^k \rightarrow 0 \:\:\: \text{as } k \rightarrow \infty, \]
	thus
	$d_H(W, C_x \times C_y) = 0$, so $W = C_x \times C_y = S'$.
\end{proof}

Equation~\eqref{eq-dimM-tangentcarpet}, Lemma~\ref{lem-carpet-tangent}  and 
Proposition~\ref{prop-tangent-dimA} combine to give
	\begin{equation*}
		\log_m(s) + \log_n(t)  = \overline{\dim}_M(W) 
			\leq \dim_A(W) 
			 = \dim_A(\hat{S}) \leq \dim_A(S).\qedhere
	\end{equation*}
\end{proof}

\section{Conformal Assouad dimension}\label{sec-cdima}
\begin{proof}[Proof of Theorem~\ref{thm-cdima}]
We first consider the case when either $s=m$, or $t=n$ (or both).
In this case we have constructed a weak tangent $W$ to $S$ that
is the product of a self-similar Cantor set and a line.
Such spaces are minimal for Assouad dimension~\cite[Lemma 6.3]{Pan-89-cdim}.
Combining Theorem~\ref{thm-main} and the second part of
Proposition~\ref{prop-tangent-dimA}, we have
\[
	\dim_A(S) = \dim_A(W) = \Cdim_A(W)
		\leq \Cdim_A(S) \leq \dim_A(S).
\]

Now we may assume that $s < m$ and $t < n$.  We wish
to show that $S$ has $\Cdim_A(S) = 0$.
As seen in~\cite[Theorem 4.1]{Tys-01-CdimA}, it suffices to show that
$S$ is \emph{uniformly disconnected}: there exists some
$C>0$ so that for every ball $B(z,r) \subset S$, there
is no $\frac{r}{C}$-chain of points joining $z$ to
$S \setminus B(z,r)$.  That is, there is no sequence 
$z=z_0, z_1, \ldots, z_N$ in $S$
so that $d(z_i, z_{i+1}) \leq \frac{r}{C}$ for
$0 \leq i < N$, with $z_N \notin B(z,r)$.

This property is not immediate since, even though $t<n$,
$S$ may project onto the unit interval in the $x$-axis.
(See Figure~\ref{fig-43dust}.)

Suppose $z \in S \cap R_k(p,q)$ (see \eqref{eq-Rk}).
Since $s < m$, any $\frac{1}{2}m^{-(k+1)}$-chain cannot travel vertically more
than $m^{-k}$.  In fact, its $y$-coordinate will stay entirely
inside either $(q/m^k,(q+2)/m^k)$ or $((q-1)/m^k,(q+1)/m^k)$.
We will assume the former, and show that suitable chains
cannot travel too far to the right.

Consider the rectangles $R_k(p+i,q)$ and $R_k(p+i,q+1)$,
for $1 \leq i \leq n$.  Since $n^{-l} \geq m^{-k}$, any
$\frac{1}{2}m^{-k}$-chain moving through these rectangles to the
right must pass through either $R_k(p+i,q)$ or $R_k(p+i,q+1)$ for
each $1 \leq i \leq n$.  As $t<n$, for some $1 \leq i \leq n$ we
must have that the interior of $R_k(p+i,q+1)$ does not meet $S$,
and so the chain passes through $R_k(p+i,q)$ from left to right.
Since $t<n$, it is impossible for any
$\frac{1}{2}n^{-(l+1)}$-chain to travel through $S \cap R_k(p+i,q)$
from left to right.

A similar argument shows that chains cannot travel too far to the left.

In summary, we have shown that $\frac{1}{2}n^{-1}m^{-k}$-chains cannot
escape from the ball $B(z,2n^{-(l-1)})$.  
(Note that $\frac{1}{2}n^{-1}m^{-k} \leq \min\{\frac{1}{2}m^{-(k+1)},\frac{1}{2}n^{-(l+1)}\}$.)
Given arbitrary $r$, we can choose $k$ so that
\[
	2n^{-(l-1)} \leq 2 n^2 m^{-k} \leq r \leq 2 n^2 m^{-(k-1)}.
\]
Therefore,
\[
	\frac{1}{2}n^{-1}m^{-k} = \frac{ 2 n^2 m^{-(k-1)}  }{4mn^3} \geq \frac{r}{4mn^3}.
\]
We have shown that for any $z \in S$, $r > 0$, no $\frac{r}{C}$-chain from $z$
can leave $B(z,r)$, where $C = 4mn^3$.
\end{proof}

\section{Lalley-Gatzouras carpets}\label{sec-lal-gatz}

As discussed in the introduction, Lalley and Gatzouras calculated the Hausdorff and upper Minkowski
dimensions of sets generalizing the construction of Bedford and McMullen.
Such a set $S$ arises as the limit set of the semigroup generated by the mappings
$A_{ij} : \R^2 \ra \R^2$ defined by
\[
	A_{ij}(x,y) = \big(a_{ij}x+c_{ij},\ b_i y + d_i \big), \quad (i,j) \in \cJ,
\]
where $\cJ = \{ (i,j) : 1 \leq i \leq m, 1 \leq j \leq n_i \}$ is the index set.
The constants are fixed to satisfy $0 < a_{ij} < b_i < 1$ for each $(i,j)$,
$\sum_{i=1}^m b_i \leq 1$, and $\sum_{j=1}^{n_i} a_{ij} \leq 1$ for each $i$.
The self-affine copies are forced to be disjoint be requiring that
$0 \leq d_1 < d_2 < \cdots < d_m < 1$ with $d_{i+1} \geq d_i + b_i$ and $1 \geq d_m+b_m$,
and, for each $i$,
$0 \leq c_{i1} < c_{i2} < \cdots < c_{in_i} < 1$ with $c_{i(j+1)} \geq c_{ij} + a_{ij}$, and $1 \geq c_{in_i}+a_{in_i}$.

Let $C_y$ be the self-similar Cantor set which is the projection of $S$ onto the $y$-axis.
Recall that its Hausdorff dimension is $\beta_y$, where $\beta_y \in (0,1]$ 
is the solution to $\sum_{i=1}^m b_i^{\beta_y}=1$.
Lalley and Gatzouras calculate the following.
\begin{theorem}[{\cite[Theorem 2.4]{Lal-Gatz-92-self-affine}}]\label{thm-lg-mink}
	The upper Minkowski dimension of $S$ is the unique $\delta$ satisfying
	$\sum_{i=1}^m \sum_{j=1}^{n_i} b_i^{\beta_y} a_{ij}^{\delta-\beta_y} = 1$.
\end{theorem}

Choose $i_* \in \{1, 2, \ldots, m\}$ so that the solution $\beta_x$ to 
$\sum_{j=1}^{n_{i_*}} {a_{i_*j}}^{\beta_x} = 1$ is maximized.
The transformations $T_j: \R \ra \R$ defined by
$T_j(x) = a_{i_*j}x+c_{i_*j}$, for $1 \leq j \leq n_{i_*}$, generate a semigroup whose limit set $C_x$ is
a self-similar Cantor set of Hausdorff dimension $\beta_x$.

The proof of Lemma~\ref{lem-carpet-tangent} easily adapts to give the following lemma.
\begin{lemma}\label{lem-carpet-tangent-lg}
	There is a weak tangent $W$ of $S$ containing a bi-Lipschitz copy of $C_x \times C_y$.
\end{lemma}
As a consequence, we have $\dim_A(S) \geq \beta_x+\beta_y$.  
Assuming Theorem~\ref{thm-main-gl} to be true, we can calculate the conformal Assouad dimension of $S$.
\begin{proof}[Proof of Theorem~\ref{thm-cdima-gl}]
	If $\beta_x = 1$ or $\beta_y = 1$, then one of $C_x$ or $C_y$ is the entire interval.
	As in the proof of Theorem~\ref{thm-cdima}, this implies that $S$ is minimal for conformal Assouad dimension.
	
	If $\beta_x < 1$ and $\beta_y < 1$, then again we can show that $S$ is uniformly disconnected,
	and so we have $\Cdim_A(S)=0$.
\end{proof}

All that remains is to complete the proof of Theorem~\ref{thm-main-gl} by
showing that $\dim_A(S) \leq \beta_x + \beta_y$.
To do this, we adapt the somewhat technical arguments of Lalley and Gatzouras 
used to prove Theorem~\ref{thm-lg-mink}, and we assume that the reader has
access to their paper.  In this proof, $C$ is a constant which varies as necessary.
\begin{proof}[Proof of Theorem~\ref{thm-main-gl}]
	First we must define the analogue of the approximate squares of \eqref{eq-Rk}.
	Note that we move between the sequence space $\cJ^\N$ and the limit set $S$ as necessary.
	
	Given $\omega = ((i_1,j_1),(i_2,j_2),\ldots) \in \cJ^\N$, and $k \in \N$,
	let $l \in \N$ be maximal so that
	\[
		R_k(\omega) := \prod_{\nu=1}^k b_{i_\nu} \leq \prod_{\nu=1}^l a_{i_\nu j_\nu}.
	\]
	Note that $l \leq k$.
	An \emph{approximate square} is the set $B_k(\omega)$ of all $\omega' = ((i_1', j_1'), \ldots) \in \cJ^\N$
	satisfying $i_\nu' = i_\nu$ for $1 \leq \nu \leq k$, and $j_\nu' = j_\nu$ for $1 \leq \nu \leq l$.
	
	As in the proof of Theorem~\ref{thm-main}, and \cite[Lemma 2.1]{Lal-Gatz-92-self-affine}, it suffices
	to show that there exists $C > 0$ so that for any $\eps >0$,
	any approximate square $B_{k'}(\omega')$ can be covered using
	at most $C(R_{k'}(\omega') / \eps)^{\beta_x+\beta_y}$
	approximate squares of diameter comparable to $\epsilon$.
	
	Following \cite[Lemma 2.2]{Lal-Gatz-92-self-affine}, it suffices to count the number of elements
	of the following set, for fixed $\omega'$, $k'$ and $l'$:
	let $\cF_\eps^*$ be the set of all
	\[
		(i_1, i_2, \ldots, i_{k+1}; j_1, j_2, \ldots, j_{l+1}),
	\]
	satisfying
	\[
		\prod_{\nu=1}^k b_{i_\nu} \geq \eps > \prod_{\nu=1}^{k+1} b_{i_\nu}
		\quad \text{and} \quad
		\prod_{\nu=1}^l a_{i_\nu j_\nu} \geq \eps > \prod_{\nu=1}^{l+1} a_{i_\nu j_\nu},
	\]
	with $i_\nu = i_\nu'$ for $\nu=1, \ldots, k'$, $j_\nu = j_\nu'$ for $\nu=1, \ldots, l'$,
	and, finally, we require one of the following two conditions to hold.
	
	\vspace{2mm}
	\noindent\textbf{Condition 1:} $1 \leq l' \leq k' \leq l+1 \leq k+1$.  Then
	\begin{enumerate}
		\item $i_\nu = i_\nu'$, $j_\nu \in \{1, \ldots, n_{i_\nu'}\}$, for $\nu = l'+1, \ldots, k'$,
		\item $(i_\nu, j_\nu) \in \cJ$, for $\nu = k'+1, \ldots, l+1$,
		\item $i_\nu \in \{1, \ldots, m\}$, for $\nu = l+2, \ldots, k+1$.
	\end{enumerate}
		
	\vspace{2mm}
	\noindent\textbf{Condition 2:} $1 \leq l' \leq l+1 \leq k' \leq k+1$.  Then
	\begin{enumerate}
		\item $i_\nu = i_\nu'$, $j_\nu \in \{1, \ldots, n_{i_\nu'}\}$, for $\nu = l'+1, \ldots, l+1$,
		\item $i_\nu = i_\nu'$, for $\nu = l+2, \ldots, k'$,
		\item $i_\nu \in \{1, \ldots, m\}$, for $\nu = k', \ldots, k+1$.
	\end{enumerate}
	
	We begin by counting the size of the subset $\cF_2$ of $\cF_\eps^*$ with Condition 2.
	Fix $R = R_{k'}(\omega')$.
	In (3), we count the set
	\[
		\Bigg\{ (i_{k'+1}, \ldots, i_{k+1}) : \prod_{\nu=k'+1}^k b_{i_\nu} \geq \frac{\eps}{R} > 
			\prod_{\nu=k'+1}^{k+1} b_{i_\nu} \Bigg\},
	\]
	which by \cite[Lemma 2.3]{Lal-Gatz-92-self-affine} has cardinality at most $C(R/\eps)^{\beta_y}$.
	The number of choices in (1) is bounded from above by a constant multiple of the number of 
	$\eps$ balls needed to cover a horizontal cross section of $S$ of length $R$,
	which is bounded from above by $C(R/\eps)^{\beta_x}$.
	These choices combine to give an upper bound of $C(R/\eps)^{\beta_x+\beta_y}$ for the size of $\cF_2$.
	
	It remains to count the size of the subset $\cF_1$ of $\cF_\eps^*$ satisfying Condition 1.
	By choosing $j_{l'+1}, \ldots, j_{k'}$, we have determined a rectangle $T$ of height $R$
	and width $R u$, where $u$ is the aspect ratio $u= \prod_{\nu = l'+1}^{k'} a_{i_\nu' j_\nu}$.
	
	The number of rectangles of width $\eps$ and height $\eps/u$ required to cover $T$
	equals the number of approximate squares of size $\eps$ needed to cover an approximate square of
	side $Ru$, which by Theorem~\ref{thm-lg-mink} is bounded by $C(Ru / \eps)^\delta \leq C (Ru/\eps)^{\beta_u+\beta_y}$.
	
	The number of approximate squares of side $\eps$ needed to cover a rectangle of width $\eps$ and
	height $\eps/u$ is at most $C(\frac{\eps/u}{\eps})^{\beta_y} = C(1/u)^{\beta_y}$,
	by \cite[Lemma 2.3]{Lal-Gatz-92-self-affine}.
		
	Combining these observations, we see that the size of $\cF_1$ is at most:
	\begin{align*}
		C\ \sum_{j_{l'+1}, \ldots, j_{k'}} \bigg( \frac{R u}{\eps} \bigg)^{\beta_x + \beta_y} 
			\left(\frac{1}{u} \right)^{\beta_y}
		& = C \left( \frac{R}{\eps} \right)^{\beta_x+\beta_y}
			\sum_{j_{l'+1}, \ldots, j_{k'}} \bigg( \prod_{\nu = l'+1}^{k'} a_{i_\nu' j_\nu} \bigg)^{\beta_x}
		\\ & \leq C \left( \frac{R}{\eps} \right)^{\beta_x+\beta_y}.
	\end{align*}
	The last inequality follows from the definition of $\beta_x$.
	
	Combining both cases, we conclude that, as desired,
	\[ | \cF_\eps^* |  \leq C \left(\frac{R}{\eps}\right)^{\beta_x+\beta_y}. \qedhere \]
\end{proof}

\bibliographystyle{plain}
\bibliography{biblio}

\end{document}